\documentclass[final,reqno]{siamltex}
\usepackage{amsmath}
\usepackage[table]{xcolor}
\usepackage{algorithm}
\usepackage{algorithmic}   
\usepackage{subfig}
\usepackage{graphicx}
\usepackage{multirow}
\usepackage{amsfonts}

\usepackage{tikz}
\usetikzlibrary{positioning, calc, decorations.pathreplacing}
\newcommand{\mybrace}[4]{ \draw[#4, decoration={brace, amplitude=6}, decorate] ([yshift=3mm]#1)--([yshift=3mm]#2) node[midway, above=2mm]{#3}; }
\usepackage[table]{xcolor}
\usepackage{booktabs}
\usepackage{epsfig,epsf,psfrag}
\usepackage[left=1in, right=1in, top=1in, bottom=1in]{geometry}

\def\x{{\bf x}}
\newcommand{\sign}[1]{\mathrm{sgn}(#1)}

\def \erf{\text{erf}}
\def \erfc{\text{erfc}}

\def\beq{\begin{equation}}
\def\eeq{\end{equation}}

\def\newpage{\vfill\eject}

\def\erf{{\rm erf}}

\def\rprime{\hbox{\hskip.25em\raise.5ex\hbox{$'$}\hskip.15em}}

\def\ddn1{{\frac{\partial}{\partial \nu_{\yb}}}}

\usepackage{geometry}
\title{ A Fast Adaptive Method for the Heat Equation with Moving or Free Boundaries in One Dimension}
\author{Chengyue Song%
\thanks{Department of Mathematics, 
Tsinghua University, Beijing, China.
Email: songcy22@mails.tsinghua.edu.cn}
\and
Jun Wang%
\thanks{Yau Mathematical Sciences Center,
Tsinghua University, Beijing, China.
Email: jwang2020@tsinghua.edu.cn.}
}

\begin{document}

\maketitle
\pagestyle{myheadings}
\markboth{\sc C. Song and J. Wang}
{\sc A Fast Adaptive Numerical Solver for Heat Equation in One Dimension}

\begin{abstract}
We present a fast adaptive method for the evaluation of heat potentials, which
plays a key role in the integral equation approach for the solution of the heat
equation, especially in a nonstationary domain. The algorithm utilizes a sum-of-exponential based fast Gauss transform that evaluates the convolution of a Gaussian with either discrete or continuous volume distributions. The latest implementation of the algorithm allows for both periodic and free space boundary conditions. The history dependence is overcome by splitting the heat potentials into a smooth history part and a singular local part. We discuss the resolution of the history part on an adaptive volume grid in detail, providing sharp estimates that
allow for the construction of an optimal grid, justifying the efficiency of the bootstrapping scheme. While the discussion in this paper is restricted to one spatial dimension, the generalization to two and three dimensions is straightforward. The performance of the algorithm is illustrated via several numerical examples.

\end{abstract}


\vspace{.1in}


	\section{Introduction}
In this paper, we discuss some key elements in the solution of the heat equation in nonstationary domain in one spatial dimension:
\begin{equation}\label{heateq}
\left\{
\begin{aligned}
\frac{\partial u}{\partial t}-\Delta u& = 0, \quad x\in \Omega(t)=[a(t),b(t)]\\
u(x,0)&=f(x),\quad x\in \Omega(t)=[a(t),b(t)]\\
u(a(t),t)&=g_{a}(t),\quad u(b(t),t) = g_b(t) \\
a(t) &< b(t), t\geq0
\end{aligned}
\right.
\end{equation}

There are many approaches for the solution of such a problem. 
Here we restrict out attention to the integral equation methods \cite{freespaceheat, lin1993thesis, wang2016njit, wang2017thesis, semilin2023}.

Classical potential theory suggests that the solution to this problem, $u$, can be written in the following form :
\begin{equation}\label{volterra}
u(x,t) = J_{\Omega(0)}[f](x,t) + D[\varphi](x,t) , \quad x \in \Omega(t)
\end{equation}
Here $J_{\Omega(0)}[f](x,t)$ is referred to as the initial heat potential and $D[\varphi](x,t)$ is referred to as the double-layer heat potential, in which $K$ is the heat kernel $\displaystyle{K(x,t)=\frac{1}{\sqrt{4\pi t}}e^{-\frac{x^2}{4t}}}$:
\begin{equation}
\begin{aligned}
J_{[a,b]}[f](x,t) &= \int_{a}^{b} K(x-y,t)f(y)dy \\
D[\varphi](x,t) &= \int_{0}^{t} \int_{\Gamma(\tau)} \frac{\partial K}{\partial \nu_y}(x-y,t-\tau)\varphi(y,\tau)ds_yd\tau
\end{aligned}
\end{equation}
It can be verified directly that any function $u(x,t)$ written in this form
readily satisfies the homogeneous heat equation and the initial condition.
To enforce the boundary condition, we let $x\rightarrow x_0\in\partial\Omega(t)$
and recall the well known jump relation of the double layer potential to obtain
an integral equation for the unknown density function $\varphi$:
\begin{equation}
\label{vol_mov}
\frac{\varphi(x,t)}{2}-D[\varphi](x,t)=-g(x,t)+J_{\Omega(0)}[f](x,t) , \quad x \in \Gamma(t)
\end{equation}

This equation is a Volterra equation of the second kind, which is well-conditioned and admits a variety of methods for its numerical solution \cite{kress, guenther1988, pogorzelski}. Once the
solution of $\varphi$ is obtained, it is substituted back into the
representation \eqref{volterra} to recover the solution to the boundary value
problem. Here we postpone the discussion of time stepping methods while observe that the fast and accurate evaluation of the heat potentials plays a key role in such numerical methods. 

The initial potential $J_{[\Omega(0)]}[f](x,t)$, when evaluated naively, costs $O(MN)$ work if $N$ denotes the number of time steps
and $M$ denotes the number of points in the spatial discretization of $\Omega(0)=[a(0),b(0)]$. Fortunately during the last few decades,
a family of fast algorithms named the fast Gauss transform (FGT) \cite{fgt1991, DFGT, strain1991vfgt, baxter2002sisc, wan2006jcp, greengard1998nfgt, beylkin2005jcp, lee2006dtfgt, tausch2009sisc, spivak2010sisc, sampath2010pfgt, veerapaneni2008jcp, wang2018sisc, fgt2023arxiv} have been developed, reducing the cost of this task to $O(N+M)$.
It is recently proposed in \cite{DFGT} that the sum-of-exponential approximation of the heat kernel readily leads to one such algorithms
which is remarkably simple and efficient especially in one dimension. In this paper, we extend the algorithm proposed in \cite{DFGT}
to the case of volume integrals with an arbitrarily small $\delta$ and incorporate periodic and free space boundary conditions.
We then demonstrate its efficacy in the evaluation of initial potentials.

When the double layer potential is considered, the task becomes even more challenging because of the full history dependence of the
time integral, costing $O(M^2N^2)$ work for direct evaluation. The FGT can be used again for the spatial integral while for the
time integral, fast methods to march in time have been proposed []. Most of these methods adopt a divide-and-conquer strategy to
decompose the layer potential $D[\varphi](x,t)$ into a local part and a history part:
\begin{equation}\label{decomposition}
\begin{aligned}
D_H[\varphi](x,t) &= \int_{0}^{t-\Delta t} \int_{\Gamma(\tau)} \frac{\partial K}{\partial \nu}(x-y,t-\tau)\varphi(y,\tau)dyd\tau \\ 
D_L[\varphi](x,t) &= \int_{t-\Delta t}^t \int_{\Gamma(\tau)} \frac{\partial K}{\partial \nu}(x-y,t-\tau)\varphi(y,\tau)dyd\tau 
\end{aligned}
\end{equation}

Here $\Delta t$ is the size of one time step. The fast algorithms for the marching in time fall into two categories: one is the Fourier
method, making use of the fact that the history part is smooth and can be represented by a few Fourier modes, and that the Fourier
modes satisfy a simple recursion in time, eliminating the history dependence problem. An obvious issue with these methods is that
they are intrinsically nonadaptive. In this paper, we adopt the framework of a fast bootstrapping method proposed in \cite{freespaceheat,lin1993thesis, wang2017thesis, semilin2023}. The key
observation is that (i) the history part $D_H$, when resolved on an adaptive grid in the volume, requires only a very sparse grid, with $O(1)$
grid points in one dimension, and $O(N_b)$ points in two and three dimensions ($O(N_b)$ being the number of grid points required
to resolve the boundary). (ii) In the spatial domain, $D_H$ also satisfies a simple recursion in time, requiring only one volume FGT and a local update in each step.
In this paper, we give a rigorous statement of property (i) and prove it in a direct fashion, producing sharp error bounds that allow
for the construction of an optimal spatial grid to resolve the heat layer potentials. A hybrid quadrature first proposed in \cite{hybquad} is
adopted for the singular quadrature in the local part. In this paper we carry out the derivation of the asymptotics to higher order.

Combining all these components mentioned above, what we present in this paper is a fast adaptive method for the efficient evaluation
of the heat potentials, which now put on a mathematically rigorous foundation and admits robust error control with no need for the
tuning of parameters. Although the discussion here is restricted to one dimension, the main idea generalizes directly to two and three
dimensions. 

	This paper is organized as follows. In section 2, we introduce simplifications of the notations for the solution of 1D heat equation and briefly review the sum-of-exponential (SOE) approximation of the heat kernel, which is the mathematical preliminaries of our work. In section 3, An SOE-based continuous FGT algorithm and its periodization are presented. In section 4, we'll discuss the properties of the solution of the second type Volterra equation \eqref{vol_mov}. 
In section 5, we deduce bounds of derivatives of the layer potentials, and discuss the adaptive grid to resolve the double layer
heat potential. Numerical examples justifying the efficacy of the algorithm will be given in section 6. Finally in section 7, we make some concluding remarks on the algorithm.
	
	\section{Mathematical Preliminaries}
	\subsection{Integral equation formulation of the heat equation in one dimension}\label{dirichlet_bc}
By direct calculation, we know that $D[\varphi]$ has the following form:
\begin{equation}\label{DLHP_mov}
\begin{aligned}
D[\varphi](x,t)& =\int_{0}^{t}\int_{\Gamma(\tau)}\frac{\partial K}{\partial \nu_y}(x-y,t-\tau) \varphi(y,\tau)ds_y d\tau\\
&=  \int_{0}^{t}\frac{-x+a(\tau)}{4\sqrt{\pi}\sqrt{(t-\tau)^{3}}}\exp\left(-\frac{(x-a(\tau))^{2}}{4(t-\tau)}\right)\varphi(a(\tau),\tau)d\tau + \int_{0}^{t}\frac{x-b(\tau)}{4\sqrt{\pi}\sqrt{(t-\tau)^{3}}}\exp\left(-\frac{(x-b(\tau))^{2}}{4(t-\tau)}\right)\varphi(b(\tau),\tau)d\tau
\end{aligned}
\end{equation}
Under the notation
$$
H(y,\tau) = \frac{y}{4\sqrt{\pi}\tau^{3/2}} \exp \left(-\frac{y^2}{4\tau}\right),\quad  \varphi_a(\tau)=\varphi(a(\tau),\tau),\quad \varphi_b(\tau)=\varphi(b(\tau),\tau)
$$
\eqref{vol_mov} can be expressed in the following form:
\begin{equation}\label{vie} 
\left\{
\begin{aligned}
&\varphi_a(t)-\int_{0}^{t}H\left(a(t)-a(\tau),t-\tau\right)\varphi_a(\tau)d\tau+\int_{0}^{t}H\left(a(t)-b(\tau),t-\tau\right)\varphi_b(\tau)d\tau=h_a(t)= -2g_a (t)+2J_{[a(0),b(0)]}[f](a(t),t)\\
&\varphi_b(t)-\int_{0}^{t}H\left(b(t)-a(\tau),t-\tau\right)\varphi_a(\tau)d\tau+\int_{0}^{t}H\left(b(t)-b(\tau),t-\tau\right)\varphi_b(\tau)d\tau=h_b(t)= -2g_b (t)+2J_{[a(0),b(0)]}[f](b(t),t)
\end{aligned}
\right.
\end{equation}

	
\subsection{SOE approximation}
The task of discrete FGT algorithm is to evaluate the following values:
$$
u_i =\sum_{j=1}^{N}e^{-\frac{\|x_i-y_j\|^2}{4t}}\alpha_j (1\leq i \leq M)
$$

Instead of the $O(MN)$ time cost of naive calculation, an FGT algorithm can evaluate them with only $O(M+N)$ time. 
There have been different algorithms of the FGT but they are all based on a certain low rank approximation of the heat kernel.

 
 

SOE approximation uses the linear combination of some exponential functions to approximate the Gaussian kernel \cite{DFGT} \cite{laplacetransform} (and notice that $K(x,t) = G(x,t)/\sqrt{4\pi t}$):
\begin{equation}\label{soe}
G(x,t) = e^{-\frac{x^2}{4t}} \approx S_n(x,t) = \sum_{k=1}^{n} w_k e^{-t_k\frac{|x|}{\sqrt{t}}}
\end{equation}

Here $\{w_k\}$ and $\{t_k\}$ are complex parameters. To compute them, we first reduce this problem into evaluating the contour integral \cite{DFGT}
$$
\frac{1}{2\pi i}\int_{\Gamma} e^z f(z) dz
$$

which can be reduced again to the approximation of $e^{z}$ on the negative real axis (\cite{talbot}), then by transforming the domain into $(-1,1)$ and then applying Carathedory-Fejer method (\cite{cf}), we can get an accurate enough rational approximation of $e^z$ on negative real axis.

After obtaining the approximation, by some easy calculation, we can recover the SOE approximation \eqref{soe}, in which the parameters $\{w_k\}$ and $\{t_k\}$ appear conjugately. To get 10-digit accuracy, one only need to set $n=12$, and even for $13$-digit accuracy (consider the complex roundoff error, this is nearly optimal), the approximation with $n=16$ can also reach it. We also notice that these parameters can be pre-computed, so that this part won't bring any extra time cost.

	\section{One-dimensional continuous FGT}
	
	\subsection{One-dimensional continuous FGT based on SOE approximation}
	\cite{DFGT} has presented a 1D discrete FGT algorithm based on SOE approximation. Just sorting the source points and the target points, a recursion formula can be easily constructed. In this paper, we'll generalize it into the continuous case, that is to say here we focus on the evaluation of
	\begin{equation}\label{FGTtask}
	u_i = J_{[a,b]}[f](x_i,t)=\frac{1}{\sqrt{4\pi t}}\int_{a}^{b} e^{-\frac{(x_i-y)^2}{4t}}f(y)dy \quad(i=1,2,\cdots,M)
	\end{equation}
	
	Here we assume that the target points $\{x_i\}$ are already sorted and can be out of the source box $[a,b]$. By SOE approximation \eqref{soe}, we have
	\begin{equation}\label{SOE_change}
	u_i \approx \frac{1}{\sqrt{4\pi t}}\sum_{k=1}^{n}w_i h_{k,i} := \frac{1}{\sqrt{4\pi t}}\sum_{k=1}^{n}w_i \int_{a}^{b} e^{-\frac{t_k}{\sqrt{t}}|x_i-y|}f(y)dy
	\end{equation}
	
	Denote $\tau_k=\dfrac{t_k}{\sqrt{t}}$, and if $a\leq x_i \leq b$, we can split $ h_{k,i}$ into 2 parts, $h_{k,i}^{+}$ and $h_{k,i}^{-}$:
\begin{equation}\label{basic2}
h_{k,i}^{+} = \int_{a}^{x_i}  e^{-\tau_k(x_i-y)}f(y)dy, \qquad h_{k,i}^{-} = \int_{x_i}^{b}  e^{-\tau_k(y-x_i)}f(y)dy
\end{equation}

Then we have the recursion formula
\begin{equation}
\begin{aligned}\label{iteration1}
h_{k,i+1}^{+} & = \int_{a}^{x_{i+1}}  e^{-\tau_k(x_{i+1}-y)}f(y)dy  = e^{-\tau_k(x_{i+1}-x_i)}h_{k,i}^{+} +  \int_{x_i}^{x_{i+1}}  e^{-\tau_k(x_{i+1}-y)}f(y)dy =: \beta_{k,i}h_{k,i}^{+} + I_{k,i}
\end{aligned}
\end{equation}

And there will analogously be
\begin{equation}\label{iteration2}
h_{k,i}^{-}=\beta_{k,i}h_{k,i+1}^{-}+\int_{x_i}^{x_{i+1}}  e^{-\tau_k(y-x_i)}f(y)dy:=\beta_{k,i}h_{k,i+1}^{-}+J_{k,i}
\end{equation}
	
When $x_i\notin [a,b]$, this recursion can still continue (and even has much simpler form). If $x_i>b$, then
\begin{equation}
h_{k,i} =  \int_{a}^{b} e^{-\tau_k (x_i-y)}f(y)dy \Rightarrow h_{k,i+1}=\beta_{k,i}h_{k,i}
\end{equation}

If $x_{i+1}<a$, similarly we have
\begin{equation}
h_{k,i}=\beta_{k,i}h_{k,i+1}
\end{equation}
	
To accelerate the process, we want to compute as few $\{I_{k,i}\}$ and $\{J_{k,i}\}$ as possible. Fortunately, Because of the conjugation of $\{t_k\}$ and $\{w_k\}$, for any $k\in \mathbb{Z}_{n}$, there will be another $k^{\prime}\in \mathbb{Z}_{n}$ such that $(w_{k^{\prime}},t_{k^{\prime}})=(\overline{w_k},\overline{t_k})$. We substitute this into the formulas above and then can easily reduce that
$$\beta_{k^{\prime},i}=\overline{\beta_{k,i}},\quad I_{k^{\prime},i}=\overline{I_{k,i}},\quad J_{k^{\prime},i}=\overline{J_{k,i}} \quad\Rightarrow\quad w_{k^{\prime}} h_{k^{\prime},i}=\overline{w_k h_{k,i}}$$
	
	Therefore, we only need to compute half of the $h_{k,i}$. Then, by doubling their real parts and discarding their imaginary parts, we can halve the run time.
	
	\subsection{FGT based on series expansion}
	
	When $t$ is very small, the error of the SOE approximation will be severely amplified due to the coefficient $1/\sqrt{4\pi t}$. In this situation, we need another way to compute \eqref{FGTtask}. In this situation, by series expansion, we can deal with this task effectively. 
	
	Since $\erf(7)<10^{-22}$ is negligible small, so for any given target point $x_0$, as long as $I=(x_0-14\sqrt{t},x+14\sqrt{t})\subset[a,b]$, we can evaluate the value by truncating the interval into $I$. If $x\in (a+14\sqrt{t},b-14\sqrt{t})$, we have (since both of the intervals can be truncated into $I$)
	\begin{equation}\label{seriesR}
	\int_{a}^{b} K(x-y,t)f(y)dy \approx \int_{\mathbb{R}} K(x-y,t)P(y)dy
	\end{equation}
	
	Here $P|_{[a,b]}$ is the piecewise polynomial approximation of $f$, and is naturally continuated onto $\mathbb{R}$, during which it won't grow super-exponentially so that it's still accurate. Then, by applying Fourier transform to the RHS of \eqref{seriesR}, we can prove that
	\begin{equation}
	u(x)=\int_{a}^{b} K(x-y,t)f(y)dy\approx  \int_{\mathbb{R}} K(x-y,t)P(y)dy = \sum_{n=0}^{\infty}\frac{d^{2n}P}{dx^{2n}}(x)\frac{t^n}{n!}
	\end{equation}
	Since $t$ is small, we just need to truncate the series into few terms and then can get very high accuracy with very fast speed.
	For the other near-boundary points (WLOG we assume $x\in (a,a+14\sqrt{t})$), we can just normalize the sub-interval into $(-1,1)$ by scaling the variable:
	\begin{equation}\label{scaledSOE}
	J_{[a,b]}[f](x,t) \approx J_{[a,a+14\sqrt{t}]}[f](x,t)  = J_{[-1,1]}\left[f\left(7\sqrt{t}(\cdot+1)+a\right)\right]\left(\frac{x-a}{7\sqrt{t}}-1,\frac{1}{49}\right)
	\end{equation}

	 $1/49$ is not a small number; so then we can call the SOE-based FGT algorithm to compute their values.
	
	\subsection{periodization}
	This FGT on a compact interval can be applied to solve Dirichlet, Neumann and Robin boundary conditions, but it isn't enough when we meet periodic boundary condition. Under this kind of circumstance, we need to periodize our FGT. Here without loss of generality, we take the example of $[a,b]=[-1,1]$ since we can analogously scale the interval from $[a,b]$ into $[-1,1]$:
\begin{equation}\label{periodizedseries}
\begin{aligned}
u_i & =\int_{\mathbb{R}} K(x_i-y;t)\tilde{f}(y)dy =\sum_{j=-\infty}^{\infty}\int_{-1}^{1} K(x-y+2j;4t)f(y-2j)dy  =\sum_{j=-\infty}^{\infty}\int_{-1}^{1} e^{-\frac{(x+2j-y)^{2}}{\delta}}f(y)dy
\end{aligned}
\end{equation}

For the positive terms, use \eqref{SOE_change} and notice that $x+2j-y \geq 2j-2 \geq 0$:
\begin{equation}\label{positivep}
\begin{aligned}
\int_{-1}^{1} e^{-\frac{(x+2j-y)^{2}}{4t}}f(y)dy &\approx \sum_{k=1}^{n}w_{k}\int_{-1}^{1} e^{-\tau_{k}\vert x+2j-y\vert}f(y)dy \\
& = \sum_{k=1}^{n}w_{k}e^{-2\tau_k j}\int_{-1}^{1} e^{-\tau_k{(x-y)}}f(y)dy \\
& = \sum_{k=1}^{n}w_{k}(e^{-2\tau_k j}h_{k,i}^{+}+e^{-2\tau_k (j-1)}h_{k,i}^{\prime +})
\end{aligned}
\end{equation}

For $j\leq -1$ terms, we similarly have
\begin{equation}\label{negativep}
\int_{-1}^{1} e^{-\frac{(x+2j-y)^{2}}{4t}}f(y)dy =  \sum_{k=1}^{n}w_{k}(e^{2\tau_k j}h_{k,i}^{-}+e^{2\tau_k (j+1)}h_{k,i}^{\prime -})
\end{equation}

here
\begin{equation}\label{newh}
h_{k,i}^{\prime +} = e^{-2\tau_k}\int_{x_i}^{1}  e^{-\tau_k(x_i-y)}f(y)dy, \qquad h_{k,i}^{\prime -} = e^{-2\tau_k}\int_{-1}^{x_i}  e^{-\tau_k(y-x_i)}f(y)dy
\end{equation}

Similar to \eqref{iteration1} and \eqref{iteration2}, we have the following recursion formulas for $h_{k,i}^{\prime +}$ and $h_{k,i}^{\prime -}$:
\begin{equation}\label{hprime_iteration}
h_{k,i+1}^{\prime +}=\beta_{k,i}h_{k,i}^{\prime +}-e^{-2\tau_k}I_{k,i},  \qquad h_{k,i}^{\prime -}=\beta_{k,i}h_{k,i+1}^{\prime -}-e^{-2\tau_k}J_{k,i}
\end{equation}

Also, 
\begin{equation}\label{hprime_initial}
\begin{aligned}
& h_{k,0}^{\prime +} =  e^{-2\tau_k}\int_{-1}^{1}  e^{-\tau_k(-1-y)}f(y)dy=\int_{-1}^{1}  e^{-\tau_k(1-y)}f(y)dy=h_{k,n+1}^{+} \\
& h_{k,n+1}^{\prime -} =  e^{-2\tau_k}\int_{-1}^{1}  e^{-\tau_k(y-1)}f(y)dy=\int_{-1}^{1}  e^{-\tau_k(y+1)}f(y)dy=h_{k,0}^{-}
\end{aligned}
\end{equation}

With \eqref{hprime_iteration} and \eqref{hprime_initial}, the values of $\{h_{k,i}^{\prime +}\}$ and $\{h_{k,i}^{\prime -}\}$ can be well evaluated. Finally, with all the results above, we can simplify the expression of \eqref{periodizedseries} into
\begin{equation}\label{truncated}
\tilde{J}[f](x,t) = \frac{1}{\sqrt{4\pi t}}\sum_{k=1}^{n}\frac{w_k}{1-e^{-2\tau_{k}}}\left( h_{k,i}^{+}+h_{k,i}^{\prime +} + h_{k,i}^{-}+h_{k,i}^{\prime -} \right)
\end{equation}

In the end of this section, we present this schematic diagram to intuitively show the data structure and the computing order in our SOE-based FGT. Here the red part is the extra computation in periodized FGT, the blue part is the extra computation for the out-of-box target points in non-periodized FGT.

\begin{figure}[h]
\centering
\begin{tikzpicture}[node distance=1.5cm]
\node at(-1.5,0){\textcolor{blue}{$\cdots$}};
\draw[-](-1,0)--(15,0);
\draw[-](0,0)--(0,0.2);
\node[below]at(0,0){\textcolor{blue}{$x_{-1}$}};
\draw[-](2,0)--(2,0.2);
\node[below]at(2,0){$x_{0}$};
\draw[-](4,0)--(4,0.2);
\node[below]at(4,0){$x_{1}$};
\draw[-](6,0)--(6,0.2);
\node[below]at(6,0){$x_{2}$};
\draw[-](8,0)--(8,0.2);
\node[below]at(8,0){$x_{n-1}$};
\draw[-](10,0)--(10,0.2);
\node[below]at(10,0){$x_{n}$};
\draw[-](12,0)--(12,0.2);
\node[below]at(12,0){$x_{n+1}$};
\draw[-](14,0)--(14,0.2);
\node[below]at(14,0){\textcolor{blue}{$x_{n+2}$}};
\node at(15.5,0){\textcolor{blue}{$\cdots$}};

\node[above]at(3,0){$\beta_0$};
\node[above]at(5,0){$\beta_1$};
\node[above]at(7,0){$\cdots$};
\node[above]at(9,0){$\beta_{n-1}$};
\node[above]at(11,0){$\beta_n$};

\node[above]at(3,1){$I_0$};
\node[above]at(5,1){$I_1$};
\node[above]at(7,1){$\cdots$};
\node[above]at(9,1){$I_{n-1}$};
\node[above]at(11,1){$I_n$};

\node[above]at(3,-1.3){$J_0$};
\node[above]at(5,-1.3){$J_1$};
\node[above]at(7,-1.3){$\cdots$};
\node[above]at(9,-1.3){$J_{n-1}$};
\node[above]at(11,-1.3){$J_n$};

\draw[dashed](2,0)--(2,2);
\draw[dashed](4,0)--(4,2);
\draw[dashed](6,0)--(6,2);
\draw[dashed](8,0)--(8,2);
\draw[dashed](10,0)--(10,2);
\draw[dashed](12,0)--(12,2);
\draw[dashed](2,-2)--(2,-0.3);
\draw[dashed](4,-2)--(4,-0.3);
\draw[dashed](6,-2)--(6,-0.3);
\draw[dashed](8,-2)--(8,-0.3);
\draw[dashed](10,-2)--(10,-0.3);
\draw[dashed](12,-2)--(12,-0.3);

\node[above]at(2,2){$h_{k,0}^{+}(=0)$};
\node[above]at(4,2){$h_{k,1}^{+}$};
\node[above]at(6,2){$h_{k,2}^{+}$};
\node[above]at(8,2){$h_{k,n-1}^{+}$};
\node[above]at(10,2){$h_{k,n}^{+}$};
\node[above]at(12,2){$h_{k,n+1}^{+}$};

\node[below]at(2,-2){$h_{k,0}^{-}$};
\node[below]at(4,-2){$h_{k,1}^{-}$};
\node[below]at(6,-2){$h_{k,2}^{-}$};
\node[below]at(8,-2){$h_{k,n-1}^{-}$};
\node[below]at(10,-2){$h_{k,n}^{-}$};
\node[below]at(12,-2){$h_{k,n+1}^{-}(=0)$};

\draw[->](3,2.4)--(3.5,2.4);
\draw[->](4.5,2.4)--(5.5,2.4);
\node[above]at(7,2.2){$\cdots$};
\draw[->](8.7,2.4)--(9.5,2.4);
\draw[->](10.5,2.4)--(11.3,2.4);

\draw[->](3.5,-2.4)--(2.5,-2.4);
\draw[->](5.5,-2.4)--(4.5,-2.4);
\node[above]at(7,-2.5){$\cdots$};
\draw[->](9.5,-2.4)--(8.7,-2.4);
\draw[->](11,-2.4)--(10.5,-2.4);

\node[above]at(2,3.5){\textcolor{red}{$h_{k,0}^{\prime +}(=h_{k,n+1}^{+})$}};
\node[above]at(4,3.5){\textcolor{red}{$h_{k,1}^{\prime +}$}};
\node[above]at(6,3.5){\textcolor{red}{$h_{k,2}^{\prime +}$}};
\node[above]at(8,3.5){\textcolor{red}{$h_{k,n-1}^{\prime +}$}};
\node[above]at(10,3.5){\textcolor{red}{$h_{k,n}^{\prime +}$}};
\node[above]at(12,3.5){\textcolor{red}{$h_{k,n+1}^{\prime +}$}};

\draw[->][red] (12,2.8)--(2,3.5);
\draw[->][red](3.2,3.9)--(3.5,3.9);
\draw[->][red](4.5,3.9)--(5.5,3.9);
\node[above]at(7,3.7){\textcolor{red}{$\cdots$}};
\draw[->][red](8.7,3.9)--(9.5,3.9);
\draw[->][red](10.5,3.9)--(11.3,3.9);

\node[below]at(2,-3.5){\textcolor{red}{$h_{k,0}^{\prime -}$}};
\node[below]at(4,-3.5){\textcolor{red}{$h_{k,1}^{\prime -}$}};
\node[below]at(6,-3.5){\textcolor{red}{$h_{k,2}^{\prime -}$}};
\node[below]at(8,-3.5){\textcolor{red}{$h_{k,n-1}^{\prime -}$}};
\node[below]at(10,-3.5){\textcolor{red}{$h_{k,n}^{\prime -}$}};
\node[below]at(12,-3.5){\textcolor{red}{$h_{k,n+1}^{\prime -}(=h_{k,0}^{-})$}};

\draw[->][red](2,-2.8)--(12,-3.5);
\draw[->][red](3.5,-3.8)--(2.5,-3.8);
\draw[->][red](5.5,-3.8)--(4.5,-3.8);
\node[below]at(7,-3.7){\textcolor{red}{$\cdots$}};
\draw[->][red](9.5,-3.8)--(8.7,-3.8);
\draw[->][red](10.7,-3.8)--(10.4,-3.8);

\node[above]at(14,0.2){\textcolor{blue}{$h_{k,n+2}$}};
\draw[->][blue](12.4,2.5)--(14,0.7);
\node[above] at(15.5,0.2){\textcolor{blue}{$\cdots$}};
\draw[->][blue](15,-0.4)--(16,-0.4);

\node[above]at(0,0.2){\textcolor{blue}{$h_{k,-1}$}};
\draw[->][blue](1.6,-2.5)--(0,-0.4);
\node[above] at(-1.5,0.2){\textcolor{blue}{$\cdots$}};
\draw[->][blue](-1,-0.4)--(-2,-0.4);

\end{tikzpicture}
\caption{The computation order of the SOE-based FGT}
\end{figure}

\section{Solving the Volterra Equation}
In the beginning of this section, we present two claims:
\begin{itemize}
\item $h$ can be exponentially expressed, which is to say that the function $p(u)=h(e^u)$ is smooth on $(0,t_0)$.
\item The function $\kappa_{ab}[\varphi](t)=\displaystyle{\int_{0}^{t}H\left(a(t)-b(\tau),t-\tau\right)\varphi(\tau)d\tau}$ can be exponentially expressed on $(0,t_0)$.
\end{itemize}

The proofs of them cover lots of redundant calculations, so we put them into the appendix.

 Back to the solution itself, If we recognize them, we can take the following hybrid form to express $\varphi$:
\begin{itemize}
\item $0<t_C$ : To initialize the solution, we let $\varphi(t)$ be a constant. (constant initialization)
\item $t_C<t<t_0$ : We use piecewise polynomials to approximate $\eta(u) = \varphi(e^u)$, here $u \in(\log t_C, \log t_0)$. (exponential form)
\item $t\geq t_0$ : We use piecewise polynomials to approximate $\varphi(t)$. (normal form)
\end{itemize}
\begin{figure}[h]
\centering
\begin{tikzpicture}[node distance=2cm]
\draw[->](0,0)--(14,0);
\draw[-](0,0)--(0,0.2);
\node[below]at(0,0){0};
\draw[-](1,0)--(1,0.2);
\node[below]at(1,0){$t_C$};
\draw[-](7,0)--(7,0.2);
\node[below]at(7,0){$t_0$};
\draw (0,0) coordinate (Origin) -- (14,0) coordinate (Ending);
    \coordinate (tc) at (1,0);
    \coordinate (t0) at (7,0);
    \mybrace{Origin}{tc}{Constant}{orange}
    \mybrace{tc}{t0}{Piecewise exponential}{red}
    \mybrace{t0}{Ending}{Piecewise normal}{blue}
\end{tikzpicture}
\caption{The expression of $\varphi$}
\end{figure}

Now based on the collocation method, we can set about solving the Volterra equation $\eqref{vol_mov}$. The constant initialization itself has asymptotically $O(\sqrt{t_C})$ error, but we deal with the later integrals in the collocation steps, the error of the integrals are $O(t_C^{3/2})$, which is acceptable as long as $t_C$ is small enough.

Without loss of generality, we consider an arbitrary subinterval $ [t_k,t_{k+1}]$. Then  denote
\begin{equation}
\varphi_a(t) = \sum_{j=0}^{M}c_{kj}^{(a)}\tilde{T}_{kj}(t), \varphi_b(t) = \sum_{j=0}^{M}c_{kj}^{(b)}\tilde{T}_{kj}(t)
\end{equation} 

Here $\tilde{T}_{kj}$ is the $j$-th order Chebyshev polynomial whose domain is scaled onto $[t_k,t_{k+1}]$. Then based on \eqref{vie}, we can construct the following linear system with respect to $\{c_{kj}^{(a)}\}$ and $\{c_{kj}^{(b)}\}$ :
\begin{equation}\label{vie_collocation}
\left\{
\begin{aligned}
&\sum_{j=0}^{M}c_{kj}^{(a)}\left(\tilde{T}_{kj}(t_{kl})+2\int_{t_k}^{t_{kl}}H\left(a(t_{kl})-a(\tau),t_{kl}-\tau\right)\tilde{T}_{kj}(\tau)d\tau\right)-2\sum_{j=0}^{M}c_{kj}^{(b)}\int_{t_k}^{t_{kl}}H\left(a(t_{kl})-b(\tau),t_{kl}-\tau\right)\tilde{T}_{kj}(\tau)d\tau\\
=&h_a(t_{kl}) -2 \int_{0}^{t_k}\left(H\left(a(t_{kl})-a(\tau),t_{kl}-\tau\right)\varphi_a(\tau)-H\left(a(t_{kl})-b(\tau),t_{kl}-\tau\right)\varphi_b(\tau)\right)d\tau \\
&2\sum_{j=0}^{M}c_{kj}^{(a)}\int_{t_k}^{t_{kl}}H\left(b(t_{kl})-a(\tau),t_{kl}-\tau\right)\tilde{T}_{kj}(\tau)d\tau+\sum_{j=0}^{M}c_{kj}^{(b)}\left(\tilde{T}_{kj}(t_{kl})-2\int_{t_k}^{t_{kl}}H\left(b(t_{kl})-b(\tau),t_{kl}-\tau\right)\tilde{T}_{kj}(\tau)d\tau\right)\\
=&h_b(t_{kl}) -2 \int_{0}^{t_k}\left(H\left(b(t_{kl})-a(\tau),t_{kl}-\tau\right)\varphi_a(\tau)-H\left(b(t_{kl})-b(\tau),t_{kl}-\tau\right)\varphi_b(\tau)\right)d\tau \\
\end{aligned}
\right.
\end{equation}

Here $\{t_{kl}\}$ are the collocation points, which are the scaled Chebyshev points of 2nd kind:
$$
t_{kl} = \frac{t_k+t_{k+1}}{2} - \frac{t_{k+1}-t_k}{2} \cos{\frac{l\pi}{L}} (0\leq l \leq L)
$$

\eqref{vie_collocation} can also be written in a matrix form:
\begin{equation}\label{vol_mat_normal}
\left({\bf{T}}+2
\begin{pmatrix}
\tilde{\bf{I}}_{aa} & -\tilde{\bf{I}}_{ab} \\
\tilde{\bf{I}}_{ba} & -\tilde{\bf{I}}_{bb}
\end{pmatrix}\right)
\begin{pmatrix}
{\bf{c}}^{(a)} \\
{\bf{c}}^{(b)}
\end{pmatrix}
= \begin{pmatrix}
{\bf{h}}^{(a)}(t_{kl}) \\
{\bf{h}}^{(b)}(t_{kl})
\end{pmatrix} +2
\begin{pmatrix}
{\bf{D}}[\varphi](a(t_{kl}),t_{k}) \\
{\bf{D}}[\varphi](b(t_{kl}),t_{k})
\end{pmatrix}
\end{equation}

Here the matrices in LHS are
$$
\tilde{\bf{I}}_{ab} =\left( \int_{t_k}^{t_{kl}}H\left(a(t_{kl})-b(\tau),t_{kl}-\tau\right)\tilde{T}_{kj}(\tau)d\tau\right)_{lj} \text{(and so do the others)}
$$

There are $2L+2$ equations ($2$ equations for each $t_{kl}$). We can solve this system $(L = M)$ or solve a least square problem $(L>M)$ to obtain the coefficients.
For the exponential part, we also split $(t_C, t_0)$ into $K=O(1)$ subintervals $\{(\tau_{k},\tau_{k+1})\} (\tau_1=t_C, \tau_{K+1}=t_0)$. Then we can construct similar linear systems:
\begin{equation}\label{vol_mat_exp}
\left({\bf{T}}+2
\begin{pmatrix}
\hat{\bf{I}}_{aa} & -\hat{\bf{I}}_{ab} \\
\hat{\bf{I}}_{ba} & -\hat{\bf{I}}_{bb}
\end{pmatrix}\right)
\begin{pmatrix}
\hat{\bf{c}}^{(a)} \\
\hat{\bf{c}}^{(b)}
\end{pmatrix}
= \begin{pmatrix}
{\bf{h}}^{(a)}(\tau_{kl}) \\
{\bf{h}}^{(b)}(\tau_{kl})
\end{pmatrix} +2
\begin{pmatrix}
{\bf{D}}[\varphi](a(\tau_{kl}),\tau_{k}) \\
{\bf{D}}[\varphi](b(\tau_{kl}),\tau_{k})
\end{pmatrix}
\end{equation}

The differences between \eqref{vol_mat_normal} and \eqref{vol_mat_exp} are the collocation points and the LHS matrices. In \eqref{vol_mat_exp}, the collocation points  $\{\tau_{kl}\}$ satisfy that
$$
\log \tau_{kl} = \frac{\log t_k+\log t_{k+1}}{2} - \frac{\log t_{k+1}- \log t_k}{2} \cos{\frac{l\pi}{L}} (0\leq l \leq L)
$$

And the matrix $\hat{\bf{I}}_{ab}$ takes the form
$$
\hat{\bf{I}}_{ab} =\left( \int_{\tau_k}^{\tau_{kl}}H\left(a(\tau_{kl})-b(\tau),\tau_{kl}-\tau\right)\hat{T}_{kj}(\tau)d\tau\right)_{lj} \text{(and so do the others)}
$$

Here $\hat{T}_{kj}$ is the $j$-th order Chebyshev polynomial whose variable is logarithmically changed and scaled:
\begin{equation}\label{expcheb}
\hat{T}_{kj}(\tau) = T_j\left(2\frac{\log\tau - \log \tau_k}{\log\tau_{k+1}-\log\tau_k}-1\right)
\end{equation}

The detailed evaluation method for the terms in \eqref{vol_mat_normal} and \eqref{vol_mat_exp} will be introduced in the next section.

\section{The evaulation of DLHP}

Solving the Volterra equation \eqref{vol_mov} is not the end; as in \eqref{volterra}, what we really need is $D[\varphi]$. Under the notation
\begin{equation}\label{I_form}
I[\gamma,\varphi](y,t) = \int_{0}^{t}H(y-\gamma(\tau),t-\tau)\varphi(\tau)d\tau = \int_{0}^{t}\frac{y-\gamma(t-\tau)}{4\sqrt{\pi}\tau^{3/2}}\exp\left(-\frac{(y-\gamma(t-\tau))^{2}}{4\tau}\right)\varphi(t-\tau)d\tau
\end{equation}

It's obvious that
\begin{equation}\label{1to2}
D[\varphi](x,t)=-I[a,\varphi_a](x,t)+I[b,\varphi_b](x,t)
\end{equation}

In the following subsections, we discuss the properties and the evaluation methods of $I[\gamma,\varphi](y)$. Without loss of generality, we assume $y>\gamma(t)$ to avoid the jump relation (there are analogous reductions and conclusions for $y<\gamma(t)$ situations).

\subsection{Hybrid quadrature method}\label{timequad}
In this section, we discuss about the evaluation of
$$
I_{a,b}[\gamma,\varphi](y) = \int_{0}^{a} H(y-\gamma(b-\tau),\tau)\varphi(b-\tau)d\tau= \int_{0}^{a} \frac{y-\gamma(b-\tau)}{4\sqrt{\pi}\tau^{3/2}}\exp \left(-\frac{(y-\gamma(b-\tau))^2}{4\tau}\right)\varphi(b-\tau)d\tau \quad(b\geq a)
$$

Here $b$ don't need to be equal to $a$, since $I_{a,b}(y)$ may be just a local correction term in the bootstrapping method introduced in later subsections.
\cite{hybquad} presents a hybrid (asymptotic+exponentially-graded mesh) method to evaluate $D[\varphi]$ in 2D situation, which also applies to our 1D situation (and with simpler form). As \cite{hybquad} suggests, we split the integral interval into two parts: $[0,\epsilon]$ and $[\epsilon,a]$, then
\begin{equation}
\begin{aligned}\label{split}
I &= \int_{0}^{a}\frac{y-\gamma(b-\tau)}{4\sqrt{\pi}\tau^{3/2}}\exp\left(-\frac{(y-\gamma(b-\tau))^{2}}{4\tau}\right)\varphi(b-\tau)d\tau \\
& = \int_{0}^{\epsilon}\frac{y-\gamma(b-\tau)}{4\sqrt{\pi}\tau^{3/2}}\exp\left(-\frac{(y-\gamma(b-\tau))^{2}}{4\tau}\right)\varphi(b-\tau)d\tau + \int_{\epsilon}^{a}\frac{y-\gamma(b-\tau)}{4\sqrt{\pi}\tau^{3/2}}\exp\left(-\frac{(y-\gamma(b-\tau))^{2}}{4\tau}\right)\varphi(b-\tau)d\tau \\
& := I_1 + I_2
 \end{aligned}
 \end{equation}

$I_1$ can be asymptotically evaluated. More precisely, with $O(\epsilon^{3/2})$ error, we have (Here $y_0 = y-\varphi(b)$)
\begin{equation}\label{asym_formula}
\begin{aligned}
I_1 &\approx \int_{0}^{2\sqrt{\epsilon}}\frac{1}{\sqrt{\pi}}\left(\frac{y_0\varphi(b)}{z^2}+\frac{\gamma^{\prime}(b)\varphi(b)-y_0\varphi^{\prime}(b)}{4}-\frac{\gamma^{\prime}(b)\varphi^{\prime}(b)z^2}{16}\right) \exp\left(-\frac{(y_0+\frac{\gamma^{\prime}(b)}{4}z^2)^2}{z^2}\right)dz\\
&= \left\{
\begin{aligned}
&\exp\left(\frac{-\alpha-|\alpha|}{2}\right)\erfc\left(\frac{|y_0|}{2\sqrt{\epsilon}}-\frac{|\gamma^{\prime}(b)|\sqrt{\epsilon}}{2}\right)\left(\frac{\sign {y_0}+\sign {\gamma^{\prime}(b)}}{4}\varphi(b)-\frac{y_0\varphi^{\prime}(b)}{4|\gamma^{\prime}(b)|}-\frac{2+|\alpha|}{4}\frac{\gamma^{\prime}(b)\varphi^{\prime}(b)}{|\gamma^{\prime}(b)|^3}\right)\\
+&\exp\left(\frac{-\alpha+|\alpha|}{2}\right)\erfc\left(\frac{|y_0|}{2\sqrt{\epsilon}}+\frac{|\gamma^{\prime}(b)|\sqrt{\epsilon}}{2}\right)\left(\frac{\sign {y_0}-\sign {\gamma^{\prime}(b)}}{4}\varphi(b)+\frac{y_0\varphi^{\prime}(b)}{4|\gamma^{\prime}(b)|}+\frac{2-|\alpha|}{4}\frac{\gamma^{\prime}(b)\varphi^{\prime}(b)}{|\gamma^{\prime}(b)|^3}\right)\\ 
+&\frac{\varphi^{\prime}(b)}{\gamma^{\prime}(b)}\sqrt{\frac{\epsilon}{\pi}}\exp \left(-\frac{y_0^2}{4\epsilon}-\frac{\alpha}{2}-\frac{\gamma^{\prime}(b)^2\epsilon}{4}\right), \quad y_0,\gamma^{\prime}(b) \neq 0, \text{ and we denote }  \alpha = y_0\gamma^{\prime}(b) \\
&\left(\frac{\gamma^{\prime}(b)\varphi(b)}{2|\gamma^{\prime}(b)|}-\frac{\gamma^{\prime}(b)\varphi^{\prime}(b)}{|\gamma^{\prime}(b)|^3}\right)\erf \left(\frac{|\gamma^{\prime}(b)|\sqrt{\epsilon}}{2}\right)+\frac{\varphi^{\prime}(b)}{\gamma^{\prime}(b)}\sqrt{\frac{\epsilon}{\pi}}\exp \left(-\frac{\gamma^{\prime}(b)^2 \epsilon}{4}\right),\quad y_0 = 0 \\
-&\frac{y_0 \varphi^{\prime}(b)\sqrt{\epsilon}}{2\sqrt{\pi}}\exp \left(-\frac{y_0^2}{4\epsilon}\right)+\left(\frac{\varphi^{\prime}(b)y_0 |y_0|}{4}+\frac{y_0 \varphi(b)}{2|y_0|}\right)\erfc\left(\frac{|y_0|}{2\sqrt{\epsilon}}\right),\quad \gamma^{\prime}(b)=0
\end{aligned}
\right.
\end{aligned}
\end{equation}

For $I_2$, we have the following two kinds of situation:
 \begin{itemize}
\item $0.02<b-a$ : Single mesh. In this situation, $I_2$ is expressed in the normal form, so that we can assume that $\varphi$ won't strongly oscillate on $(b-a,b-\epsilon)$. Then, by the deduction in \cite{hybquad}, it's practicable to change the variable $\tau = e^{-u}$ and obtain an exponentially-graded mesh for $I_2$. 

\item $b-a\leq0.02$ : Dual mesh. In this situation, when $\tau>b-0.02$,  exponential form will appear in the expression of $I_2$, which means we have to simultaneously deal with the strong oscillation for both $H(y,\tau)$ (when $\tau \rightarrow \epsilon$) and $\varphi(b-\tau)$(when $\tau\rightarrow b$). To individually deal with the singularity on both sides of the interval, one can split $(\epsilon,a)$ into 2 sub-intervals: $(\epsilon,c)$ and $(c,a)$.
$$
I_2 = I_{21} + I_{22} = \int_{\epsilon}^{c}\frac{y}{4\sqrt{\pi}\tau^{3/2}}\exp\left(-\frac{y^{2}}{4\tau}\right)\varphi(b-\tau)d\tau + \int_{c}^{a}\frac{y}{4\sqrt{\pi}\tau^{3/2}}\exp\left(-\frac{y^{2}}{4\tau}\right)\varphi(b-\tau)d\tau
$$

The singularity of $H$ on $(\epsilon,c)$ in $I_{21}$ can be dealt with exponentially-graded mesh $\tau = e^{-u}$. From \eqref{expcheb}, the singularity of $\varphi$ on $(c,a)$ in $I_{22}$ can be well solved with another type of exponentially-graded mesh: $b-\tau =e^{-u}$. 
 
In some occasions, $b-a$ can be very small (especially on the first time step, where $b=a=t_1$), and we'll meet the constant-initialization part of $\varphi$. But this won't be a tough dilemma, since we just need to furtherly split $(c,a)$ into $(c,b-t_C)\cup(b-t_C,a)$. Constant approximation for $\varphi$ on $\tau\in(b-t_C,a)$ is extremely coarse; but since $a - (b-t_C)\leq t_C$, as long as $t_C$ is small, the error of the quadrature process of the last term will be negligibly small. So the actual expression of $I_2$ in this situation is
 \begin{equation}\label{expmesh_comp_2}
 \begin{aligned}
I_2(y,t) & = \int_{-\log c}^{-\log \epsilon} \frac{y e^{u/2}}{4\sqrt{\pi}}\exp \left(-\frac{y^2 e^{u}}{4}\right)\varphi(b-e^{-u})du  \quad\text{(singularity of $H$)}\\
& + \int_{-\log (b-c)}^{-\log t_C} \frac{y}{4\sqrt{\pi}e^{u}(b-e^{-u})^{3/2}}\exp \left( -\frac{y^2}{4(b-e^{-u})} \right) \varphi (e^{-u})du \quad\text{(singularity of $\varphi$)}\\
& + \int_{b-t_C}^{a}\frac{y}{4\sqrt{\pi}\tau^{3/2}}\exp(-\frac{y^{2}}{4\tau})\varphi(b-\tau)d\tau \quad\text{(initialization part)}
\end{aligned}
 \end{equation}
\end{itemize}

Now consider the choice of $c$. According to Sec.4, we need to deal with the exponentially oscillation (in other words, some kind of ``singularity'') of $\varphi$ on and only on $(0,t_0)$. So, when the interval is long enough, we can split the 2nd mesh ($I_{22}, b-\tau = e^u$) so that $b-\tau\in(0,t_0)$, which is to say, $c = b-t_0$. But when it's short, we shouldn't forget the singularity of  $H(y,\cdot)$ when $\tau\rightarrow 0$. In this kind of situation, we let $c = (a+\epsilon)/2$ to balance the singularity of both parts when $\tau$ is close to $\epsilon$ or $a$. The following figure presents a more intuitive introduction:

\newpage

\begin{figure}[h]
\centering
 \begin{tikzpicture}[scale=40]
		\draw[->](0,0)--(0.11,0)node[below]{$a$};
		\draw[->](0,0)--(0,0.11)node[left]{$b$};
		\node[below left](0,0){$O$};
		
		\filldraw[lightgray](0.08,0.1)--(0.003,0.023)--(0.043,0.043)--(0.1,0.1);
	\filldraw[gray](0.003,0.003)--(0.003,0.023)--(0.043,0.043);
	\filldraw[darkgray](0.003,0.023)--(0.003,0.1)--(0.08,0.1);
	\draw[fill=gray,domain=0:0.1]plot(\x,\x);
		\draw[domain=0:0.077]plot(\x+0.003,\x+0.023);
		\draw[fill=gray,domain=0.003:0.043]plot(\x,{\x/2+0.0215});
		\draw[fill = gray,densely dashed](0.003,0.003)--(0.003,0.1);
		\node[below]at(0.1,0.08){$b = a$};
		\node[below]at(0.072,0.116){$b = a+t_0$};
		\node[below]at(0.063,0.043){$b = t_0+\dfrac{a+\epsilon}{2}$};
\end{tikzpicture}
\end{figure}

\begin{center}
\fcolorbox{black}{white}{\textcolor{white}{111111}} = asymptotic only/impossible , 
\fcolorbox{black}{darkgray}{\textcolor{darkgray}{111111}} = single mesh,

\fcolorbox{black}{gray}{\textcolor{gray}{111111}} = dual mesh, $c = (a+\epsilon)/2$,
\fcolorbox{black}{lightgray}{\textcolor{lightgray}{111111}} = dual mesh, $c = b - t_0$
\end{center}

\subsection{Time stepping}

The DLHP function has time dependence. When the time variable $t$ grows, the quadrature process for $D[\varphi]$ has to set more nodes. For example, for equispaced temporal steps $t_k = k\Delta t$, If we directly evaluate the values $\{D[\varphi](y,t_k)\}$ on $N$ time points without any fast algorithms, there will be unavoidable $O(N^2)$ cost.
To accelerate this process, we'll use the bootstrapping method. More exactly, we first decompose $D[\varphi](y,t)$ into 2 parts - a history part and a local correction part (for the notation, see \eqref{decomposition}).
\begin{equation}
D[\varphi] = D_L[\varphi]+D_H[\varphi]
\end{equation}

Generally we have
\begin{equation}\label{time_stepping}
\begin{aligned}
D_H[\varphi](x,t)=&\int_0^{t-\Delta t} \int_{\Gamma(\tau)} \frac{\partial K}{\partial \nu_y}(x-y, t-\tau) \varphi(y, \tau) d s_{y} d \tau \\
=&\int_0^{t-\Delta t} \int_{\Gamma(\tau)}\left(\int_{\mathbb{R}} K(x-z, \Delta t) \frac{\partial K}{\partial \nu_y}(z-y, t-\Delta t-\tau) dz\right) \varphi(y, \tau) d s_{y} d \tau \\
=& \int_{\mathbb{R}} K(x-z, \Delta t)\left(\int_0^{t-\Delta t } \int_{\Gamma(\tau)}  \frac{\partial K}{\partial \nu_y}(z-y, t-\Delta t -\tau) \varphi(y, \tau) d s_{y} d \tau\right) d z\\
=& \int_{\mathbb{R}} G(x-z, \Delta t)D[\varphi](z,t-\Delta t ) d z
\end{aligned}
\end{equation}

This equation takes the form of the free-space FGT. On the other side, by integration by parts,
\begin{equation}\label{I_diff0}
\begin{aligned}
I[\gamma,\varphi](y,t) &=\int_{0}^{t} \left(\frac{\partial}{\partial \tau}\frac{1}{2}\erf\left(\frac{y-\gamma(\tau)}{2\sqrt{t-\tau}}\right)+K(y-\gamma(\tau),t-\tau)\gamma^{\prime}(\tau)\right)\varphi(\tau)d\tau \\
& = \frac{1}{2}\left(\varphi(t)-\erf\left(\frac{y-\gamma(0)}{2\sqrt{t}}\right)\varphi(0) \right)  -\frac{1}{2} \int_{0}^{t} \erf\left(\frac{y-\gamma(\tau)}{2\sqrt{t-\tau}}\right)\varphi^{\prime}(\tau)d\tau + \int_{0}^{t}K(y-\gamma(\tau),t-\tau)\gamma^{\prime}(\tau)\varphi(\tau)d\tau
\end{aligned}
\end{equation}

Observe that  its every term in \eqref{I_diff0} (and so that their sum $|I[\gamma,\varphi](y,t)|$ and furtherly $|D[\varphi](z,t)|$) decays super-exponentially with respect to $y$ when it gets far away from the boundary. So, we can assume that $D$ is compactly supported on $[z_L,z_R]=[a(t)-14\sqrt{t},b(t)+14\sqrt{t}]$, then we can truncate the interval from $\mathbb{R}$ into $[z_L,z_R]$ and then call the standard continuous FGT introduced above. As for the local term $D_L[\varphi]$, the hybrid quadrature method introduced in \ref{timequad} works well on it.

Also, we emphasize that these methods will be applied not only in the evaluation of $D[\varphi]$, but also in solving the Volterra equations \eqref{vol_mat_normal} and \eqref{vol_mat_exp}. For example, the elements in the LHS matrices can be evaluated by hybrid quadrature; and the RHS historical terms is DLHP function, which meets our fast time stepping algorithm above.

\subsection{Construct the spatial mesh}
\eqref{time_stepping} infers that we have to spatially approximate the DLHP function $z\mapsto D[\varphi](z,t-\Delta t)$ for the FGT. To do this, a spatial mesh is needed to evaluate the function values on them.
However, a too simple spatial mesh like an equispaced one isn't a good idea. $I(y)$ decays very quickly when $|y|\ll 0$ and approximately equals to $0$ when $|y|\gg 0$ (see later in this subsection). If we want to obtain piecewise approximation for $I[\varphi](y)$, then we have to construct a better and more adaptive mesh. 

Here our method is : 
\begin{itemize}
\item When $t$ is small, every sub-interval can be $O(\sqrt{t})$ long.
\item When $t$ is big, every sub-interval can be  $O(\sqrt[n+1]{\ln t})$ (hence practically $O(1)$) long.
\end{itemize}

Based on the reductions of the exponential expression in the appendix, here we present a direct proof, which can be generalized to higher-dimensional cases too. 

For any arbitrary sub-interval $(y_a,y_b)$, note that $I\in C^{\infty}(y_a,y_b)$ as long as $y_a>\gamma(t)$. So, consider its Taylor series about $y_0 = (y_a+y_b)/2$, when truncated into the $n-$th term (which means $I$ is approximated with a polynomial of degree $n$), the value of the Taylor reminder can be bounded as
$$
\varepsilon_n(y) = \frac{|I^{(n+1)}(\xi)|}{(n+1)!}|y-y_0|^{n+1} (\xi \in (y_a,y_b))
$$

When $t$ is small, based on the inequality \eqref{Idiffs_bound} in the appendix, denote $\alpha= (b-a)/ \sqrt{t}$, $\varepsilon_n$ can be bounded:
\begin{equation}\label{errdiff}
\begin{aligned}
|\varepsilon_{n}(y)|&\leq \frac{\max_{(a,b)}|I^{(n+1)}(\xi)|}{(n+1)!}\left(\frac{\alpha \sqrt{t}}{2}\right)^{n+1}\\
&\leq\frac{1}{(n+1)!}\left(\frac{\alpha\sqrt{t}}{2}\right)^{n+1} \left(\sum_{k\geq1}^{n}|A_k|\frac{C}{\sqrt{\pi}}\frac{2^{k/2}\sqrt{k!}}{(2\sqrt{t})^{k+1}} \exp\left(-\frac{(y-\gamma(0))^2}{8t}\right)+\lambda+\mu\sqrt{t}\right) \\
& =  \frac{\alpha^{n+1}C}{(n+1)!\sqrt{\pi}}\sum_{k=1}^{n}|A_k|\frac{\sqrt{k!}\ t^{(n-k)/2}}{2^{n+k/2+2}}\exp\left(-\frac{(y-\gamma(0))^2}{8t}\right)+\frac{\lambda+\mu \sqrt{t}}{(n+1)!}\left(\frac{\alpha\sqrt{t}}{2}\right)^{n+1}
\end{aligned}
\end{equation}

This is a polynomial with respect of $\sqrt{t}$. Since $t$ is small, this is a rather tight bound so that the approximation on $O(\sqrt{t})$ length sub-interval will be accurate enough.

However, when $t$ gets larger, the error bound in \eqref{errdiff} will blow up, and this method may be not accurate enough. In this situation, we try another way, and here our solution is to densitify the mesh from the original $O(\sqrt{t})$ one into the $O(1)$ one. Back to \eqref{Idiffs}, the bound of $S_1$ is still valid, since it decays when $t$ grow larger. But \eqref{I_bound},\eqref{S_bound} and \eqref{Idiffs_bound} present an $O(\sqrt{t})$ bound, which is not so practical. Fortunately, there is another way for estimating the bound:

For $S_2$, we first consider $I^{\prime}[\gamma,\psi]$ with $|\psi|\leq M$ on $[0,t]$. Then by splitting the interval into $(0,t-\epsilon)\cup(t-\epsilon,t)$, we have
\begin{equation}\label{S2_part1}
\begin{aligned}
|I_{1}^{\prime}| &= \left|\frac{\partial}{\partial y}\int_{0}^{\epsilon}\frac{y-\gamma(t-\tau)}{4\sqrt{\pi}\tau^{3/2}}\exp\left(-\frac{(y-\gamma(t-\tau))^2}{4\tau}\right)\psi(t-\tau)d\tau\right| 
\end{aligned}
\end{equation}

The derived function can be asymptotically evaluated by \eqref{asym_formula}. Under our assumption $y_0>0$, so it's differentiable. On the other hand, in each situation every term and its derivative in \eqref{asym_formula} are $O(1)$, so $|I_1|$  is $O(1)$ too.
\begin{equation}\label{S2_part2}
\begin{aligned}
|I_2^{\prime}| &= \left|\int_{\epsilon}^{t}\frac{\partial}{\partial y}\left(\frac{y-\gamma(t-\tau)}{4\sqrt{\pi}\tau^{3/2}}\exp\left(-\frac{(y-\gamma(t-\tau))^2}{4\tau}\right)\right)\psi(t-\tau)d\tau\right| \\
\leq & M\int_{\epsilon}^{t}\left|1-\frac{(y-\gamma(t-\tau))^2}{2\tau}\right|\frac{1}{4\sqrt{\pi}\tau^{3/2}}\exp\left(-\frac{(y-\gamma(t-\tau))^2}{4\tau}\right)d\tau\\
\leq & M\int_{\epsilon}^{t}\left(1+\frac{2}{e}\right)\frac{1}{4\sqrt{\pi}\tau^{3/2}}d\tau <\left(1+\frac{2}{e}\right)\frac{1}{2\sqrt{\pi\epsilon}}
\end{aligned}
\end{equation}

From \eqref{S2_part1} and \eqref{S2_part2}, it's obvious that : no matter how large $t$ is, $|S_2|=O(1)$. 

For $S_3$, we similarly discuss $I[\gamma,\psi]$ with $|\psi|\leq M$ on $[0,t]$:
\begin{equation}
|I_1| = \left|\int_{0}^{\epsilon}\frac{y-\gamma(t-\tau)}{4\sqrt{\pi}\tau^{3/2}}\exp\left(-\frac{(y-\gamma(t-\tau))^2}{4\tau}\right)\psi(t-\tau)d\tau\right|=(O(1)) \quad\left(\text{by \eqref{asym_formula}}\right)
\end{equation}
\begin{equation}
\begin{aligned}
|I_2| =& \left|\int_{\epsilon}^{t}\frac{y-\gamma(t-\tau)}{4\sqrt{\pi}\tau^{3/2}}\exp\left(-\frac{(y-\gamma(t-\tau))^2}{4\tau}\right)\psi(t-\tau)d\tau\right| \\
\leq & \int_{\epsilon}^{t}\frac{1}{\sqrt{8\pi e}}\frac{1}{\tau} Md\tau = \frac{M(\ln t-\ln \epsilon)}{\sqrt{8\pi e}}
\end{aligned}
\end{equation}

So, $S_3=O(\ln t)$. Analogous to \eqref{errdiff}, it's easy to show that an $O(\sqrt[n+1]{\ln t})$ mesh is enough to bound $|\varepsilon_n(y)|$ into $O(1)$ level. From a practical point of view, $O(\sqrt[n+1]{\ln t})$ grows so slowly when $t$ grows larger so that it has almost no difference to $O(1)$ on any arbitrary time interval $[0,T]$. For example, when we take $n=16$, then$\sqrt[17]{\ln 10^{8}} \approx 1.187$, which means we just need only a little (or even no) densification to the spatial mesh.

In the end, in both situations, we notice the bound decays super-exponentially with respect to $y$. Therefore, when $y$ (more precisely, $y/\sqrt{t}$) grows larger, the error bound of each order will be lower, so that the approximation will be more precise, which means the mesh can be coarser (i.e. larger $\alpha$) when it's far from $y=0$.

Based on all the results above, we can look back to the time cost of the bootstrapping method of $N$ steps:
\begin{itemize}
\item When time steps forward, the spatial mesh will be won't be denser. Hence, less and less spatial nodes are needed to approximate the DLHP function. From this perspective, there are only $O(1)$ source and target points (the spatial nodes for the piecewise approximation of $D[\varphi]$) in each run of the free-space FGT, which means the time cost for FGT is $O(1)$ each time.
\item The evaluation of every local correction term obviously takes $O(1)$ time.
\end{itemize}

So , we can claim that the bootstrapping method manages to accelerate the time stepping process from $O(N^2)$ time to $O(N)$ time.

\section{Numerical examples}
\subsection{Numerical results for the FGT}
	In the following examples, we choose the source function $f(y)=\sin(10\pi y)$. To test the periodic and nonperiodic FGT, by choosing different $n$ and $t$ and fix $10^{6}$ equispaced points on $[-1,1]$ as the target points, we have following results:
	
	\begin{minipage}{\textwidth}

\begin{minipage}[t]{0.48\textwidth}
\makeatletter\def\@captype{table}
\begin{tabular}{cccc}
\toprule
$n$\qquad&$t$\qquad&$L_{\infty}$ Error\qquad&Time cost (s) \\
\midrule
8&$1$&$1.6\times10^{-8}$&0.68 \\
12&$1$&$3.7\times10^{-12}$&0.92 \\
16&$1$&$2.2\times10^{-12}$&1.03\\
\\
8&$10^{-1}$&$2.4\times10^{-11}$&0.67 \\
12&$10^{-1}$&$4.0\times10^{-13}$&0.85 \\
16&$10^{-1}$&$1.8\times10^{-11}$&1.02\\
\\
8&$10^{-2}$&$2.3\times10^{-8}$&0.67 \\
12&$10^{-2}$&$1.2\times10^{-12}$&0.75 \\
16&$10^{-2}$&$2.0\times10^{-12}$&0.93\\
\\
8&$10^{-3}$&$1.7\times10^{-8}$&0.57 \\
12&$10^{-3}$&$1.7\times10^{-12}$&0.75 \\
16&$10^{-3}$&$1.0\times10^{-11}$&0.93 \\
\\
8&$10^{-3}$&$1.5\times10^{-9}$&0.56 \\
12&$10^{-3}$&$1.3\times10^{-12}$&0.75 \\
16&$10^{-4}$&$2.1\times10^{-12}$&0.93 \\
\bottomrule
\end{tabular}
\caption{Periodic case}
\end{minipage}
\begin{minipage}[t]{0.48\textwidth}
\makeatletter\def\@captype{table}
\begin{tabular}{cccc}
\toprule
$n$\qquad&$t$\qquad&$L_{\infty}$ Error\qquad&Time cost (s) \\
\midrule
8&$1$&$1.6\times10^{-8}$&0.49 \\
12&$1$&$3.8\times10^{-12}$&0.64 \\
16&$1$&$7.4\times10^{-12}$&0.78\\
\\
8&$10^{-1}$&$2.6\times10^{-9}$&0.48 \\
12&$10^{-1}$&$1.3\times10^{-12}$&0.63 \\
16&$10^{-1}$&$2.1\times10^{-11}$&0.77\\
\\
8&$10^{-2}$&$2.0\times10^{-8}$&0.48 \\
12&$10^{-2}$&$9.7\times10^{-13}$&0.64 \\
16&$10^{-2}$&$2.9\times10^{-12}$&0.77\\
\\
8&$10^{-3}$&$1.6\times10^{-8}$&0.48 \\
12&$10^{-3}$&$1.8\times10^{-12}$&0.64 \\
16&$10^{-3}$&$1.0\times10^{-11}$&0.77 \\
\\
8&$10^{-4}$&$2.7\times10^{-9}$&0.49 \\
12&$10^{-4}$&$1.3\times10^{-12}$&0.64 \\
16&$10^{-4}$&$2.1\times10^{-12}$&0.78 \\
\bottomrule
\end{tabular}
\caption{Nonperiodic case}
\end{minipage}
\end{minipage}

\subsection{Numerical results for the Volterra equation solver module}

Now we test the accuracy when solving the Volterra equation \eqref{vie}. To do this, we consider the following options for \eqref{vie} and solve it on $\Omega_T = (a(t),b(t))\times(0,1)$ with the method introduced in section 4:
\begin{equation}
\left\{
\begin{aligned}
a(t) &= -1+\text{J}_2(10t)\\
b(t) &= \frac{\sin (20t)}{3}\\
h_a(t) &=\sin(10\pi t)\\
h_b(t) &= 1-\text{J}_1({10t})
\end{aligned}
\right. ,\quad t_0 = 0.02,\quad L = M+1
\end{equation}
here $J_\nu$ is the Bessel function of the first kind.

 From the integral equation theory we know that when the collocation order $L$ is fixed, there will be a $k>0$ so that the error $\varepsilon$ and the time step length $h $ have the following relationship: $\varepsilon = O(h^{k})$ . By choosing different collocation orders and different time steps, we have the following results:
\begin{table}[h]
\centering
\begin{tabular}{ccccccc}
\toprule
\multicolumn{3}{c}{Exponential}& &
		\multicolumn{3}{c}{Normal}\\
		\midrule
Order&$k$&$R^2$&&Order&$k$&$R^2$\cr
		\cmidrule{1-3}\cmidrule{5-7}                 
		4&-3.40 &0.9969&&6 &-6.09 &0.9621\\
		6&-5.63 &0.9988&&8 &-7.70 &0.9604\\
		8&-8.76 &0.9983&&10 &-9.88 &0.9706\\
\bottomrule
\end{tabular}
\caption{The relationship between the error $\varepsilon_{M,N}$ and the collocation order $M$, the number of the time steps $N$}
\end{table}

\subsection{Numerical results for the periodic heat equation solver}
We consider the periodic problem example
\begin{equation}\label{periodic_eq}
\left\{
\begin{aligned}
u_t-u_{xx}& = F(x,t)=\left(-k\pi \sin (k\pi t)+(k+1)^2\pi^2\cos(k\pi t)\right)\sin (k+1)\pi x , \\
u(x,0)&=\sin (k+1)\pi x \\
\end{aligned}
\right.
\end{equation}

Here $k\in\mathbb{N}$, and $(x,t)\in[-1,1]\times[0,1]$. To solve this problem, we first choose $n$ equispaced time steps $t_j=j/N (1\leq j\leq N)$.
For the periodic problem, we know
$$
u(x,t)=J[\tilde{f}]+\tilde{V}[F]=\int_{\mathbb{R}}K(x-y,t)\tilde{f}(y)dy+\int_{0}^{t}\int_{\mathbb{R}}K(x-y,t-\tau)\tilde{F}(y,\tau)dyd\tau\quad(\tilde{f} and \tilde{F}\text{ are periodized onto } \mathbb{R} \text{ by } f \text{ and } F)
$$
and there's the local correction formula
\begin{equation}\label{periodic_solution}
u(x, t+\Delta t)=\int_{\mathbb{R}} u(z, t) K(x-z, \Delta t) d z +\int_t^{t+\Delta t} \int_{\mathbb{R}} K(x-y, t+\Delta t-\tau) \tilde{F}(y, \tau) d y d \tau
\end{equation}

To compute \eqref{periodic_solution}, we can compute the first term by (periodic) FGT, and then discretize the second term with respect to $\tau$ on every step $(t_i,t_{i+1})$ (such as Gauss-Legendre quadrature). Then, we reduce the periodic heat equation into some periodic FGTs. For different parameter $k$, by taking different time step length and different discretization order, the numerical result approximately shows that $\epsilon=O(A^{-n})$ . More precisely, we have the following results:
\begin{table}[h]
\centering
\begin{tabular}{ccccccccccc}
\toprule
\multicolumn{3}{c}{$k=6$}& &
		\multicolumn{3}{c}{$k=8$}& &\multicolumn{3}{c}{$k=10$}\\
		\midrule
$L$&$A$&$R^2$&&$L$&$A$&$R^2$&&$L$&$A$&$R^2$\cr
		\cmidrule{1-3} \cmidrule{5-7} \cmidrule{9-11}
		4&1.126 &0.9988& &4 &1.070 &0.9988& &4 &1.042 &0.9987\\ 
		8&1.624 &0.9990& &8 &1.341 &0.9995& &8 &1.206 &0.9991\\ 
		16&6.679 &0.9992& &16 &3.349 &0.9999& &16 &2.229 &0.9998\\ 
\bottomrule
\end{tabular}
\caption{The relationship between the error $\epsilon$ and the quadrature discretization order $L$, the number of time steps $n$}
\vspace{-0.5cm}
\end{table}

To quantify the error, we compare the $u(\cdot,1)$ solved by the solver and its exact solution $u(x,t)=\cos k\pi t\cdot\sin(k+1)\pi t$ and calculate the $L_2$ error between them to obtain $\epsilon$.

\subsection{Numerical results for the nonperiodic heat equation solver}
For the nonperiodic case, we consider the Dirichlet problem on $(a(t),b(t))\times(0,0.5)$:
\begin{equation}\label{nonperiodic_eq}
\left\{
\begin{aligned}
u_t-\Delta u& = 0, \quad x\in \Omega(t)=(a(t),b(t))\\
u(x,0)&=f(x)=\sin(k\pi x),\quad x\in \Omega(t)\\
u(a(t),t)&=\int_{-2}^{2}\frac{1}{\sqrt{4\pi t}}\exp \left(-\frac{(a(t)-y)^2}{4t}\right)\sin (k\pi y) dy \\
u(b(t),t)&=\int_{-2}^{2}\frac{1}{\sqrt{4\pi t}}\exp \left(-\frac{(b(t)-y)^2}{4t}\right)\sin (k\pi y) dy \\
a(t) &= \frac{\sin(k\pi t)}{2},\quad b(t) = 1-\log(1+t) + J_1(k\pi t)
\end{aligned}
\right.
\end{equation}
It's easy to show that the following integral is the solution of \eqref{nonperiodic_eq} :
\begin{equation}\label{nonperiodic_solution}
u(x,t)=\int_{-2}^{2}\frac{1}{\sqrt{4\pi t}}\exp \left(-\frac{(x-y)^2}{4t}\right)\sin (k\pi y) dy, \quad (x,t)\in (a(t),b(t))\times(0,0.2)
\end{equation}
We'll respectively test the exponential part and the normal part, to show the relationship between the error and the time steps:
\begin{itemize}
\item for the exponential part, we take $t_0=0.02$ and compare the $u(\cdot,t_0)$ solved by the solver and its exact solution \eqref{nonperiodic_solution}, and then calculate the $L_2$ error between them to obtain $\epsilon$. For different parameter $k$ and different collocation order $L$, the numerical result approximately shows that $\epsilon=O(K^{-p})$, here $K$ is the number of the pieces of the exponential mesh.
\begin{table}[h]
\centering
\begin{tabular}{ccccccccccc}
\toprule
\multicolumn{3}{c}{$k=6$}& &
		\multicolumn{3}{c}{$k=8$}& &\multicolumn{3}{c}{$k=10$}\\
		\midrule
$L$&$p$&$R^2$&&$L$&$p$&$R^2$&&$L$&$p$&$R^2$\cr
		\cmidrule{1-3} \cmidrule{5-7} \cmidrule{9-11}
		8& 8.22 &0.9529& &8 & 7.89 &0.9731& &8 & 6.83 &0.9749\\ 
		16& 10.78 &0.9819& &16 &10.96 &0.9803& &16 &11.18 &0.9706\\ 
		32& 11.72 &0.9911& &32 &11.91 &0.9953& &32 &12.02 &0.9689\\ 
\bottomrule
\end{tabular}
\caption{The relationship between the error $\epsilon$ and the collocation order $L$, the number of (exponential) time steps $K$}
\vspace{-0.5cm}
\end{table}
\item for the normal part, we still take $t_0=0.02$, and now consider $u(\cdot,0.5)$. We still focus on $L_2$ error here. Also, to exclude the error of the exponential part, we'll choose high-order collocation and efficient time steps on it (so that in this example, the error are almost all from the normal part). Similarly, it can be seen that $\epsilon=O(N^{-q})$, here $N$ is the number of the time steps of normal mesh. Here are the numerical results:
\begin{table}[h]
\centering
\begin{tabular}{ccccccccccc}
\toprule
\multicolumn{3}{c}{$k=6$}& &
		\multicolumn{3}{c}{$k=8$}& &\multicolumn{3}{c}{$k=10$}\\
		\midrule
$L$&$q$&$R^2$&&$L$&$q$&$R^2$&&$L$&$q$&$R^2$\cr
		\cmidrule{1-3} \cmidrule{5-7} \cmidrule{9-11}
		4& 3.80 &0.9813& &4 & 4.16 &0.9472& &4 & 3.31 &0.9338\\ 
		8& 8.27 &0.9788& &8 &9.52 &0.9857& &8 &7.97 &0.9796\\  
\bottomrule
\end{tabular}
\caption{The relationship between the error $\epsilon$ and the collocation order $L$, the number of (exponential) time steps $K$}
\vspace{-0.5cm}
\end{table}
\end{itemize}

\subsection{A physical application : solving Stefan problem}
The Stefan problem is a simple but meaningful mathematical model of phase change. In one-dimensional case, it can be represented as the melting process of a semi-inﬁnite solid on $[s(t),\infty)$, during which we assume that there isn't any volume change
occurring. At the ﬁxed boundary $x = 0$, there's a heat reservoir providing the flux for the latent heat need for melting. Then, we have the following PDE system to mathematically describe the process:
\begin{equation}\label{stefan}
\left\{
\begin{aligned}
u_t(x,t) & =u_{x x}(x, t), \quad x \in(0, s(t)) \\
u(0, t) & =T_c\left(t\right)>T_m,\  u(s(t), t) =T_m, \\
u(x,0) &= f(x),\ s(0)=s_0, \\
s^{\prime}(t) & =-\beta u_x(s(t), t)
\end{aligned}
\right.
\end{equation}
\begin{figure}[h]
\centering
\begin{tikzpicture}[node distance=2cm]
\draw[->](0,0)--(14,0);
\draw[-](0,0)--(0,0.2);
\node[below]at(0,0){$u(0,t)=T_c(t)$};
\node[below]at(3.7,0){The temperature is $u(x,t)$};
\draw[-](7,0)--(7,0.2);
\node[below]at(7,0){$s(t)$};
\draw (0,0) coordinate (Origin) -- (14,0) coordinate (Ending);
    \coordinate (s) at (7,0);
    \mybrace{Origin}{s}{Liquid phase}{red}
    \mybrace{s}{Ending}{Solid phase}{blue}
    \node[below]at(10.75,0){The temperature  is $T_m$};
\draw[purple][->](6.3,0.6)--(7.7,0.6);
\node[above]at(7,0.6){\textcolor{purple}{melting}};
\end{tikzpicture}
\caption{Schematic of one-dimensional one-phase Stefan problem}
\end{figure}
\vspace{0cm}

To solve $u(x,t)$ on $t=t_k (k=1,2,\cdots,n)$, we can combine the heat equation solver constructed above with the SDC method \cite{sdc} to solve the ODE with respect to $s$:
\begin{algorithm}[h]
	\caption{Solve the Stefan problem \eqref{stefan} by SDC method}
	\begin{algorithmic}[1]
	\FOR {$i = 1,2,3,\cdots,n$}
		\STATE Generate Chebyshev nodes $\{t_{ij}\}$ on $(t_{i-1},t_i)$ (here $t_0=0$) and then use forward-Euler discretization to obtain an initial approximation for $\{u(x,t_{ij})\}$ and $\{s(t_{ij})\}$
		\STATE Use SDC method to obtain more accurate solution for $\{u(x,t_{ij})\}$
		\STATE Interpolate to obtain $s(t)|_{(t_{i-1},t_i)}$
		\ENDFOR
		\RETURN $s(t)$ (defined on $(0,t_n)$) and $\{u(x,t_k)\}(k=1,2,\cdots n)$
	\end{algorithmic}  
\end{algorithm}

We first consider the following Stefan problem on $(0,s(t))\times(0,1)$:
\begin{equation}\label{stefan_t0}
\left\{
\begin{aligned}
u_t(x,t) & =u_{x x}(x, t), \quad x \in(0, s(t)) \\
u(0, t) & =u_0,\  u(s(t), t) =0,  \\
s(0)&=2\lambda\sqrt{t_0}, \\
u(x,0) &= u_0 \left(1-\frac{\erf\left(\frac{x}{2\sqrt{t_0}}\right)}{\erf(\lambda)}\right)\\
s^{\prime}(t) & =-\beta u_x(s(t), t)
\end{aligned}
\right.
\end{equation}

Here $\lambda$ satisfies that $\lambda e^{\lambda^2}\erf(\lambda)=\frac{u_0}{\beta \sqrt{\pi}}$. This problem has an explicit solution
\begin{equation}\label{stefan_explicit}
\left\{
\begin{aligned}
u(x,t) &= u_0 \left(1-\frac{\erf\left(\frac{x}{2\sqrt{t+t_0}}\right)}{\erf(\lambda)}\right) \\
s(t)& = 2\lambda\sqrt{t+t_0}
\end{aligned}
\right.
\end{equation}

As an example, by taking the parameters as $\lambda = 1/2, \beta=1, t_0=0.1$, our method gives numerical results in the following table. Here we quantify the accuracy of the method by considering the error of $s(1)$ (denote it as $\varepsilon_{n,l}$). Then we can similarly observe that the relationship $\varepsilon_{n,l}=O(h^{-k})$ approximately holds, here $h$ is the step length:
\begin{table}[h]
\centering
\begin{tabular}{ccccccccccc}
\toprule
\multicolumn{3}{c}{$l=4$}& &
		\multicolumn{3}{c}{$l=7$}& &\multicolumn{3}{c}{$l=10$}\\
		\midrule
Order&$k$&$R^2$&&Order&$k$&$R^2$&&Order&$k$&$R^2$\cr
		\cmidrule{1-3} \cmidrule{5-7} \cmidrule{9-11}
		4&3.01 &0.9999& &4 &3.49 &0.9994& &4 &3.51 &0.9995\\
		8&2.65 &0.9999& &8 &3.52 &0.9988& &8 &5.89 &0.9991\\
		12&2.61 &0.9999& &12 &3.65 &0.9996& &12 &4.78 &0.9998\\
\bottomrule
\end{tabular}
\caption{The relationship between the error $\varepsilon_{n,l}$ and the Chebyshev discretization order $n$, the SDC iteration round $l$}
\end{table}

Then, fixing these parameters, we modify the flux $T_c$ into $T_c^{\prime}=u_0\cdot (1+\cos(5\pi t))/2$ to consider a much more complicated condition, in which the modified problem has no more explicit solution. However, our method is still feasible. By choosing different time steps, we have the following results (the notations are the same as above):

\newpage

\begin{table}[h]
\centering
\begin{tabular}{ccccccccccc}
\toprule
\multicolumn{3}{c}{$l=4$}& &
		\multicolumn{3}{c}{$l=8$}& &\multicolumn{3}{c}{$l=12$}\\
		\midrule
Order&$k$&$R^2$&&Order&$k$&$R^2$&&Order&$k$&$R^2$\cr
		\cmidrule{1-3} \cmidrule{5-7} \cmidrule{9-11}
		8& 1.53 & 0.8972 & &8 & 4.50 &0.9952 & &8 &4.43 &0.9960\\
		12& 1.52 & 0.8989& &12 & 3.31& 0.9786& &12 &3.70 &0.9738\\
		16& 1.53 & 0.9025& &16 & 2.62& 0.9860& &16 &4.23 &0.9738\\
\bottomrule
\end{tabular}
\caption{The relationship between the error $\varepsilon_{n,l}$ and the Chebyshev discretization order $n$, the SDC iteration round $l$}
\end{table}

When $l=4$, the SDC iteration doesn't converge yet, so the $R^2$ coefficient is not so high. But when $l$ get larger, the SDC iteration converges, and there holds that $\varepsilon_{n,l}=O(h^{-k})$.

\begin{figure}[h]
\centering
\includegraphics[scale=0.5]{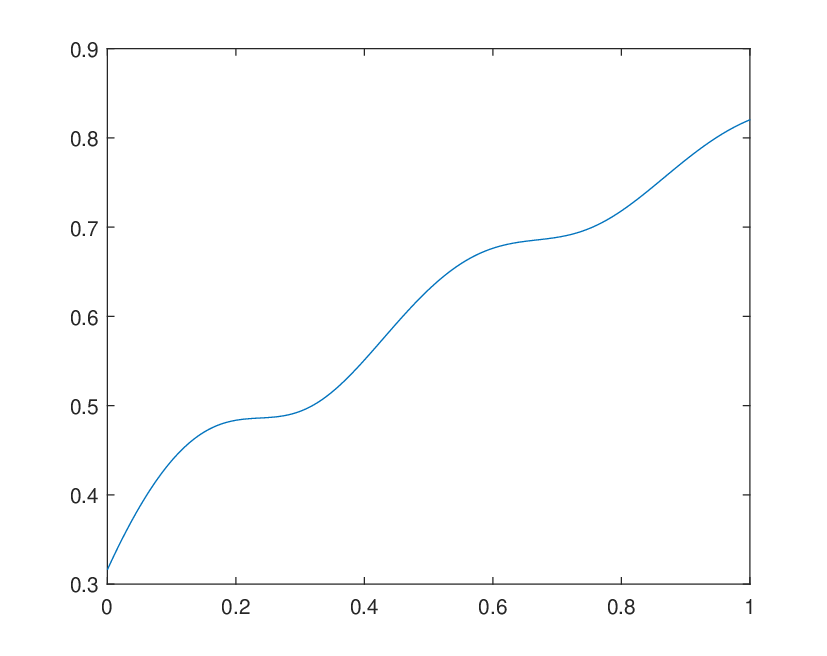}
\caption{The figure of solution $s(t)$ of the modified Stefan problem example}
\end{figure}

\section{Conclusion}
The fast and accurate evaluation of heat layer potentials is an essential ingredient in the construction of efficient integral equation
methods for the solution of the heat equation. In this paper, we have examined two main aspects of this task.
An SOE based fast Gauss transform for discrete and volume sources is presented, and a fast bootstrapping method for the marching
in time is discussed in detail.
Furtherly, by combining the algorithm with a proper method for the numerical solution of the Volterra equation, we present an efficient method
for the moving boundary problems of the heat equation.

\newpage


\appendix
\section{The exponential expressions}
In Sec.4, we present two claims, based on which it can be readliy seen that the solution $\varphi$ to the Volterra equations \eqref{vol_mov} can be exponentially expressed. We give their proofs here.
\subsection{The exponential expression for $h$}
From \eqref {vie}, without loss of generality, we only need to consider $h_a(t)=J_{\Omega (0)}[f](a(t),t)-g_{a}(t)$, and then get similar conclusion can be given for $h_b = J_{\Omega (0)}[f](b(t),t)-g_{b}(t)$.
For $g_a(t)$, when we consider $\eta_{g}(u)=g_a(e^u)$,  we have
\begin{equation}\label{gexp}
\eta_{g_a}^{(k)}(u) = \sum_{j=1}^{k}e^{ju}g_a^{(j)}(e^{u}) = \sum_{j=1}^{k}t^{j}g_a^{(j)}(t) (t<t_0)
\end{equation}

So, due to the fact that $t_0 \ll 1$, the terms in the RHS decays exponentially with respect to $j$, which means each order of the derivatives of $\eta_a$ can be well bounded by different orders of the deratives of $g_a$ (and so that $\eta_b$ and $g_b$).

For $J_a(t)=J_{\Omega (0)}[f](a(t),t)$, the situations will be more complicated. First we define $J(x,t) = J_{[a,b]}[f](x,t)$, obviously $\displaystyle{\frac{\partial J[f]}{\partial t}=\frac{\partial ^{2}J[f]}{\partial x^{2}}}$, so by integration by parts, we have
\begin{equation}
\begin{aligned}\label{h_sumform}
\frac{\partial^{n} J}{\partial x^{n}}(x,t) &= \sum_{k=1}^{n-1} (-1)^{n-k}\frac{\partial^{k}K}{\partial x^{k}}(x-y,t)\bigg|_{y=a}^{y=b} + (-1)^{n}J[f^{(n)}](x,t)
\end{aligned}
\end{equation}

By Cramer's inequality, we have (here $w = z/2\sqrt{t}$):
\begin{equation}\label{cramer}
\left|\frac{\partial^{k}K}{\partial y^{k}}(z,t)\right| = \frac{1}{\sqrt{\pi}}\frac{1}{(2\sqrt{t})^{k+1}}\left|\frac{d^{k}}{dw^{k}}(e^{-w^2})\right| \leq \frac{C}{\sqrt{\pi}}\frac{2^{k/2}\sqrt{k!}}{(2\sqrt{t})^{k+1}} \exp\left(-\frac{z^2}{8t}\right)
\end{equation}

Where $C\leq1.09$. Then substitute this into \eqref{h_sumform}, and then
\begin{equation}\label{h_bound2}
\begin{aligned}
\left|\frac{\partial^{n} J}{\partial x^{n}}(x,t)\right| &\leq\sum_{k=1}^{n-1} \frac{C}{\sqrt{\pi}}\frac{2^{k/2}\sqrt{k!}}{(2\sqrt{t})^{k+1}}\left( \exp\left(-\frac{(x-a)^2}{8t}\right)+\exp\left(-\frac{(x-b)^2}{8t}\right)\right) + \left|J[f^{(n)}](x,t)\right|
\end{aligned}
\end{equation}

Since
 \begin{equation}\label{total_diff2}
 \begin{aligned}
\frac{\text{d}^nJ_{a}}{\text{d}t^n} &= \sum_{p+q\leq n}J^{(p,q)}(a(t),t)\left(\sum_{\substack{\sum i\alpha_i =n-q \\ \sum \alpha_i = p}} \frac{n!}{q!} \prod_{i\geq 1} \frac{a^{(i)}(t)^{\alpha_i}}{(i \alpha_i)!}\right) \\
& =\sum_{p+q\leq n}\frac{\partial^{p+2q} J}{\partial x^{p+2q}} (a(t),t)\left(\sum_{\substack{\sum i\alpha_i =n-q \\ \sum \alpha_i = p}} \frac{n!}{q!} \prod_{i\geq 1} \frac{a^{(i)}(t)^{\alpha_i}}{(i \alpha_i)!}\right) := \sum_{p+q\leq n}S_{n,p,q}[a]\frac{\partial^{p+2q} J}{\partial x^{p+2q}} (a(t),t)
\end{aligned}
\end{equation}

Finally, similar to \eqref{gexp}, for $J_a$ we have
\begin{equation}\label{Jexp}
\begin{aligned}
\left|\eta_{J_a}^{(k)}(u)\right|& = \sum_{j=1}^{k}t^{j}\cdot\left|\frac{\text{d}^{j}J_{a}}{\text{d}t^j}(a(t),t)\right|  =   \sum_{j=1}^{k}t^{j}\sum_{p+q\leq j}\left|S_{n,p,q}[a]\right|\left|\frac{\partial^{p+2q} J}{\partial x^{p+2q}}(a(t),t)\right| \\
& \leq  \sum_{j=1}^{k}t^{j}\sum_{p+q\leq j}|S_{n,p,q}[a]|\left(\sum_{l=1}^{p+2q-1} \frac{C}{\sqrt{\pi}}\frac{2^{l/2}\sqrt{l!}}{(2\sqrt{t})^{l+1}}\left( e^{-\frac{(a(t)-a)^{2}}{8t}}+e^{-\frac{(a(t)-b)^{2}}{8t}}\right) + \left|J[f^{(n)}](x,t)\right|\right)
\end{aligned}
\end{equation}

Consider the power of $t$ in each term of the RHS. It's easy to see the power of $t$ in every term is $j-\dfrac{l+1}{2}$, and since $l\leq p+2q-1, p+q\leq j$, we have
\begin{equation}\label{power_bound1}
j-\frac{l+1}{2} \geq j-\frac{p}{2}-q \geq  \frac{p}{2}
\end{equation}
So we've presented a good bound for $\left|\eta_{J_a}^{(k)}(u)\right|$ when $t<t_0\ll 1$. 

From \eqref{gexp} and $\eqref{Jexp}$, by changing the variable, under the notation $\eta_a(u) = \eta_{J_a}(u)-\eta_{g_a}(u)$, each order of the derivatives of $\eta_a$ is well bounded. Of course, these two formulas also infer that when $t$ grow larger (for example, $\geq 1$), such expressions can't reduce the derivatives any more. So in this situation, we'll change back to normal expression.

\subsection{The exponential expression for $\kappa$}
For \eqref{I_form}, it's clear that  by integration by parts we have
\begin{equation}
\begin{aligned}\label{I_diff1}
\frac{\partial I}{\partial y}[\gamma,\varphi](y,t)&=\int_{0}^{t}\left(\frac{\partial}{\partial \tau}K(y-\gamma(\tau),t-\tau)-\gamma^{\prime}(\tau)\frac{y-\gamma(\tau)}{4\sqrt{\pi}(t-\tau)^{3/2}}\exp\left(-\frac{(y-\gamma(\tau))^2}{4(t-\tau)}\right)\right)\varphi(\tau)d\tau \\
& = -K(y-\gamma(0),t)\varphi(0)-\int_{0}^{t}K(y-\gamma(\tau),t-\tau)\varphi^{\prime}(\tau)d\tau-I[\gamma,\gamma^{\prime} \varphi](y)
\end{aligned}
\end{equation}

When $n\geq 2$, based on the trick used in \eqref{h_sumform}, there will be the following recursion formula:
\begin{equation}\label{recursion}
\frac{\partial^{n}I}{\partial y^{n}}[\gamma,\varphi](y)=- \frac{\partial^{n-1}}{\partial y^{n-1}} K(y-\gamma(0),t)\varphi(0)+ \frac{\partial^{n-2}I}{\partial y^{n-2}}[\gamma,\varphi^{\prime}](y) -\frac{\partial^{n-1}I}{\partial y^{n-1}}[\gamma,\gamma^{\prime}\varphi](y)  
\end{equation}

By \eqref{recursion} and some combinatorial reductions, It's clear that $I^{(n)}$ can be expressed as the linear combination of these three kinds of terms:
\begin{equation}\label{Idiffs}
\begin{aligned}
\frac{\partial^{n}I}{\partial y^{n}}[\gamma,\varphi](y)=&\sum_{k\geq1}^{n-1}\frac{\partial^{k}K}{\partial y^{k}} (y-\gamma(0),t)\sum_{\substack{\sum (i-1)\alpha_i =p-r \\ \sum \alpha_i = q \\ 2p+q = n-k}}(-1)^{q+1}(p+q)!\prod_{i\geq 1}\frac{\gamma^{(i)}(0)^{\alpha_i}}{(i \alpha_i)!}\varphi^{(r)}(0) \quad (:=S_1) \\
+& \sum_{\substack{\sum (i-1)\alpha_i =p-r \\ \sum \alpha_i = q \\ 2p+q = n-1}}\frac{(-1)^{q}(p+q)!}{\prod_{i\geq 1}(i \alpha_i)!}\frac{\partial I}{\partial y}\left[\gamma,\prod_{i=1}^{m}\gamma^{(i)\alpha_i}\varphi^{(r)}\right](y) \quad (:=S_2) \\
+& \sum_{\substack{\sum (i-1)\alpha_i =p-r \\ \sum \alpha_i = q \\ 2p+q = n}}\frac{(-1)^{q}(p+q)!}{\prod_{i\geq 1}(i \alpha_i)!}I\left[\gamma,\prod_{i=1}^{m}\gamma^{(i)\alpha_i}\varphi^{(r)}\right](y) \quad (:=S_3)
\end{aligned}
\end{equation}

It seems complicated, but the terms in \eqref{Idiffs} can be classified into 3 types: the linear combination of the spatial deratives of $K$ ($S_1$),  the linear combination of the deratives of some DLHP functions ($S_2$), and the linear combination of some DLHP functions ($S_3$).

For $S_1$, Cramer's inequality \eqref{cramer} gives an efficient bound (here for simplicity, we denote the coefficient for the $k-$th term of $S_1$ as $A_k$):
\begin{equation}
S_1 \leq\sum_{k\geq1}^{n-1}|A_k|\frac{C}{\sqrt{\pi}}\frac{2^{k/2}\sqrt{k!}}{(2\sqrt{t})^{k+1}} \exp\left(-\frac{(y-\gamma(0))^2}{8t}\right)
\end{equation}

One another thing we need to claim is that $\varphi^{(r)}(0)=0$  since we choose constant for initialization (so that $S_1$ won't blow up as long as the inputs $f$ and $g$ are smooth enough). 

For $S_2$, \eqref{I_diff1} infers that this kind of terms can be expressed in the form of the linear combination of heat kernels, SLHP functions and DLHP functions. 

For $S_3$ (and the DLHP functions transformed from $S_2$), we have the following bound: 
\begin{equation}\label{I_bound}
\begin{aligned}
|I[\gamma,\varphi](y,t)| &\leq \frac{|\varphi(0)|+|\varphi(t)|}{2}+\frac{1}{2}\int_{0}^{t}|\varphi^{\prime}(\tau)|d\tau + |S[\gamma^{\prime}\varphi](y,t)|
\end{aligned}
\end{equation}

Here $S$ is the single-layer heat potential function. Then for any $S[\varphi]$ we also have the following bound
\begin{equation}\label{S_bound}
|S[\varphi](z,t)| = \left|\int_{0}^{t}K(z-\gamma(\tau),t-\tau)\varphi(\tau)d\tau\right| \leq\int_{0}^{t}\frac{1}{\sqrt{4\pi (t-\tau)}}|\varphi(\tau)|d\tau \leq \sqrt{\frac{t}{\pi}}\max_{\tau\in[0,t]} |\varphi(\tau)|
\end{equation}

\eqref{I_bound} and \eqref{S_bound} means that $|S_2|$ and $|S_3|$ are both $O(1+\sqrt{t})$ functions, which means there exists two constants $\lambda,\mu>0$ so that
\begin{equation}\label{Idiffs_bound}
\begin{aligned}
\left|\frac{\partial^{n}I}{\partial y^{n}}[\gamma,\varphi](y)\right|\leq&\sum_{k\geq1}^{n-1}|A_k|\frac{C}{\sqrt{\pi}}\frac{2^{k/2}\sqrt{k!}}{(2\sqrt{t})^{k+1}} \exp\left(-\frac{(y-\gamma(0))^2}{8t}\right) + \lambda + \mu \sqrt{t}
\end{aligned}
\end{equation}

Finally, under our notation we have $\kappa_{aa}(t) = I[a,\varphi](a(t),t)$. Under the change of the variable $\eta_{\kappa}(u)=\kappa(e^u)$, by imitating the analysis in \eqref{gexp},\eqref{total_diff2} and \eqref{Jexp}, there will be the following bound:
\begin{equation}\label{kappaexp}
\begin{aligned}
\left|\eta_{\kappa}^{(k)}(u)\right|& = \sum_{j=1}^{k}t^{j}\cdot\left|\frac{\text{d}^{j}\kappa_{aa}}{\text{d}t^j}(a(t),t)\right| \\
& =   \sum_{j=1}^{k}t^{j}\sum_{p+q\leq j}\left|S_{n,p,q}[y]\right|\left|\frac{\partial^{p+2q} I}{\partial y^{p+2q}}[a,\varphi](a(t),t)\right| \\
& \leq  \sum_{j=1}^{k}t^{j}\sum_{p+q\leq j}|S_{n,p,q}[y]|\left(\sum_{k\geq1}^{p+2q-1}|A_k|\frac{C}{\sqrt{\pi}}\frac{2^{k/2}\sqrt{k!}}{(2\sqrt{t})^{k+1}} \exp\left(-\frac{(a(t)-a(0))^2}{8t}\right) + \lambda + \mu\sqrt{t}\right)
\end{aligned}
\end{equation}

Finally, we consider the power of $t$ in every term of  \eqref{kappaexp}. Similar to \eqref{power_bound1}, since
\begin{equation}
j -\frac{p+2q}{2}\geq p+q-\frac{p+2q}{2}\geq \frac{p}{2}
\end{equation}

We can say that \eqref{kappaexp} also gives a good bound for $\left|\eta_{\kappa}^{(k)}(u)\right|$ when $t<t_0\ll 1$. 

 So as we've claimed in the beginning of this section, piecewise polynomials can be used to approximate $\eta(u) = \varphi(e^u)$, here $u < \log t_0$.

\newpage
\bibliographystyle{siam}
\bibliography{ref.bib}
\end{document}